\documentclass[a4paper,12pt]{article}
\usepackage{amsmath,amssymb,fullpage,mathptmx}
\DeclareMathOperator{\KS}{\mathit{C}}
\DeclareMathOperator{\KP}{\mathit{K}}

\begin{document}
\author{Andrej A. Muchnik\thanks{The game interpretation considered in this paper was suggested by Andrej Muchnik (24.02.1958--18.03.2007) in his talks at Kolmogorov seminar (Moscow Lomonosov Univerisity). The examples are added and text is prepared by I.~Mezhirov (\texttt{mezhirov@gmail.com}, University of Kaiserslautern), A.~Shen (\texttt{alexander.shen@lif.univ-mrs.fr}, LIF Marseille, CNRS \& University Aix--Marseille, on leave from IITP, RAS, Moscow) and N.~Vereshchagin (\texttt{nikolay.vereshchagin@gmail.com}, Moscow State Lomonosov University) who are responsible for all errors and omissions. Preparation of this paper was supported in part by ANR Sycomore and  NAFIT ANR-08-EMER-008-0[1,2] grants and RFBR 09-01-00709-Á grant.}, \ Ilya~Mezhirov, Alexander Shen, Nikolay Vereshchagin}
\title{Game interpretation of Kolmogorov complexity}
\date{}
\maketitle
\begin{abstract}
The Kolmogorov complexity function $K$ can be relativized using any oracle $A$, and most properties of $K$ remain true for relativized versions  $K^A$. In section~\ref{game-interpretation} we provide an explanation for this observation by giving a game-theoretic interpretation and showing that all ``natural'' properties are either true for all $K^A$ or false for all $K^A$  if we restrict ourselves to sufficiently powerful oracles~$A$. This result is a simple consequence of Martin's determinacy theorem, but its proof is instructive: it shows how one can prove statements about Kolmogorov complexity by constructing a special game and a winning strategy in this game.
\end{abstract}

\section{Game interpretation}
\label{game-interpretation}

Consider all functions defined on the set of binary strings and having  non-negative integer values, i.e., the set $\mathcal{F}=\mathbb{N}^{\{0,1\}^*}$. Let $\alpha$ be a property of such a function (i.e., a subset of $\mathcal{F}$). We say that $\alpha$ is $O(1)$-\emph{stable} if $f_1\in\alpha\Leftrightarrow f_2\in\alpha$ for any two functions $f_1,f_2\in \mathcal{F}$ such that $f_1(x)=f_2(x)+O(1)$, i.e., the difference $|f_1(x)-f_2(x)|$ is bounded.

Let $A$ be an oracle (a set of strings). By $K^A(x)$ we denote the Kolmogorov complexity of a string $x$ relativized to oracle $A$, i.e., the length of the shortest description for $x$ if the decompressor is allowed to use $A$ as an oracle. (See~\cite{li-vitanyi} or \cite{shen-uppsala} for more details; we may use either plain complexity (denoted usually by $C$ or $\textit{KS}$) or prefix complexity (denoted usually by $K$ or $\textit{KP}$) though the game interpretation would be slightly different; see below.)

For a given $A$ the function $K^A$ is defined up to  an $O(1)$ additive term, therefore an $O(1)$-stable property $\alpha$ is well defined for $K^A$ (does not depend on the specific version of $K^A$). So $\alpha(K^A)$ becomes a property of the oracle~$A$. It may be true for some oracles and false for other ones. For example, if $w_n$ is a $n$-bit prefix of Chaitin's random real $\Omega$, the ($O(1)$-stable) property ``$K^A(w_n)>0.5n+O(1)$'' is true for trivial oracle $A=\mathbf{0}$ and false for $A=\mathbf{0}'$. The following result (a special case of a classical result of D.~Martin~\cite{martin-0}, see the discussion below) shows, however,  that for ``usual'' $\alpha$ the property $\alpha(K^A)$ is either true for all sufficiently large $A$ or false for all sufficiently large~$A$.

\medskip

\textbf{Proposition}. \emph{Let $\alpha$ be a Borel property. Then there exists an oracle $A_0$ such that either
$\alpha(K^A)$ is true for all $A\ge_T A_0$
or
$\alpha(K^A)$ is false for all $A\ge_T A_0$.}

\medskip

Here $\ge_T$ stands for Turing reducibility. The statement is true for different versions of complexity (plain complexity, prefix complexity, decision complexity, a priori complexity, monotone complexity etc.). We provide the proof for plain complexity $\KS$ and there describe the changes needed for other versions.

\medskip

\textbf{Proof}. Consider the following infinite game with full information. Two players called (as usual) Alice and Bob enumerate graphs of two functions $A$ and $B$ respectively; arguments and values of $A$ and $B$ are binary strings. The players' moves alternate; at each move player may add finitely many pairs to the graph of her/his function but cannot delete the pairs that are already there (so the values of $A$ and $B$ that are already defined remain unchanged).

The winner is declared as follows. Let $K_A$ and $K_B$ be the complexity functions that correspond to decompressors $A$ and $B$, i.e.,
     $$
K_A(x)=\min\{ l(p) \mid A(p)=x\}
     $$
where $l(p)$ stands for the length of~$p$; the function $K_B$ is defined in a similar way.
Let us agree that Alice wins if the function
	$$
K(x)=\min (K_A(x), K_B(x))
        $$
satisfies~$\alpha$. If not, Bob wins. (A technical correction: functions $K_A$ and $K_B$ may have infinite values; we assume that $\alpha$ is somehow extended to such functions, e.g., is false for all functions with infinite values.)

\smallskip

\textbf{Lemma}. \emph{If Alice has a computable winning strategy in this game,
then $\alpha(\KS)$ is true for \textup(plain\textup) complexity function $\KS$; if Bob has a computable winning strategy, then $\alpha(\KS)$ is false.}

\smallskip

\textbf{Proof of the Lemma} is straightforward. Assume that Alice has a computable winning strategy. Let her use this strategy against the enumeration of the graph of optimal decompressor function (so $K_B(x)=\KS(x)$ for all~$x$). Note that in fact Bob ignores the moves of Alice and enumerates the graph of $B$ at its own pace. Since both playes use computable strategies, the game is computable. Therefore $K_A\le K_B+O(1)$ due to the optimality of $B$, and
    $$\min(K_A(x),K_B(x))=K_B(x)+O(1)=\KS(x)+O(1).$$
Since Alice wins and $\alpha$ is $O(1)$-stable, the function $\KS$ has property $\alpha$. The same argument (with exchanged roles of Alice and Bob) can be used if Bob has a winning strategy.~\raisebox{-0.7ex}{\hbox{\small$\Box$}}

\smallskip
The statement and the proof of the lemma can be relativized: if Alice/Bob has a winning strategy that is $A$-computable for some oracle $A$, then $\alpha(\KS^A)$ is true/false.

Now recall Martin's theorem on the determinacy of Borel games: the winning condition of the game described is a Borel set (since $\alpha$ has this property), so either Alice or Bob has a winning strategy in the game. So if the oracle~$A$ is powerful enough (is above the strategy in the hiearchy of $T$-degrees), the property $\alpha(K^A)$ is true (if Alice has a winning $A$-computable strategy) or false (if Bob has a winning $A$-computable strategy). Theorem is proven.~\raisebox{-0.7ex}{\hbox{\small$\Box$}}

\section{Discussion}
\label{discussion}

Let us make several remarks.

\medskip
\textbullet\quad
As we have said, this proposition is a consequence of an old general result proved by Martin. Lemma on p.~688 of~\cite{martin-0} together with Borel determinacy~\cite{martin-1,determinacy} guarantees that for every Borel Turing-invariant property $\Phi$ of infinite binary sequences either $\Phi$ is true for all sequences in some upper cone (in the degrees semilattice), or $\Phi$ is false for all sequences in some upper cone. It remains to note that the property $\varphi(A)=\alpha(K^A)$ is a Turing-invariant Borel property.

The proof in~\cite{martin-0} uses a different (and simpler) game: two players alternate adding bits to a sequence, and the referee checks whether the resulting infinite sequence satisfies $\Phi$. The advantage of our game is that it is more tailored to the definition of Kolmogorov complexity and therefore can be used as a prototype of games needed to prove some specific statements about Kolmogorov complexity.

\medskip
\textbullet\quad
Note that not all theorems in algorithmic information theory are $O(1)$-stable. For example, most of the results about algorithmic properties of complexity function are not stable. (The non-computablity of the complexity function or its upper semicomputablity is not a stable property, while the non-existence of a nontrivial computable lower bound is stable. Also the Turing-completeness of $\KS$ is a non-stable assertion though the stronger claim ``any function that is $O(1)$-close to $\KS$ can be used as an oracle to decide halting problem'' is stable.) The other assumption (Borel property) seems less restrictive: it is hard to imagine a theorem about Kolmogorov complexity where the property in question won't be a Borel one by construction.

\medskip
\textbullet\quad
One may ask whether the statement of our theorem can be used as a practical tool to prove the properties of Kolmogorov complexity. The answer is yes and no at the same time.  Indeed, it is convenient to use some kind of game while proving results about Kolmogorov complexity, and usually the argument goes in the same way: we let the winning strategy play against the ``default'' strategy of the opponent and the fact that the winning strategy wins implies the statement in question. However, it is convenient to consider more special games. For example, proving the inequality
        $$
\KS (x,y) \ge \KS(x)+\KS(y|x) - O (\log n)
       $$
(for strings $x$ and $y$ of length at most $n$), we would consider a game where Alice wins if $K_B(x,y)<k+l$ implies that either $K_A(x)<k+O(\log n)$ or $K_A(y|x)<l+O(\log n)$ for every $n,k,l$ and for all strings $x,y$ of length at most~$n$.

This example motivates the following version of the main theorem. Let $\alpha$ be a property of two functions in $\mathcal{F}$, i.e., a subset of $\mathcal{F}\times\mathcal{F}$. Assume that $\alpha$ is monotone in the following sense: if $\alpha (f,g)$ is true, $f'(x)\le f(x)+O(1)$, and $g'(x)\ge g(x)-O(1)$, then $\alpha(f',g')$ is true, too. Consider the version of the game when Alice wins if $\alpha(K_A,K_B)$ is true. If Alice has a computable winning strategy, then $\alpha(\KS,\KS)$ is true; if Bob has a computable winning strategy, then $\alpha(\KS,\KS)$ is false. (The proof remains essentially the same.)

\smallskip
We provide several examples where game interpretation is used to prove statements about Kolmogorov complexity in Section~\ref{sec:examples};
other examples can be found in~\cite{noncompressible} and in the survey~\cite{ver-survey}.

\medskip
\textbullet\quad
Going in the other direction, one would like to extend this result to arbitrary results of computablility theory not necessarily related to Kolmogorov complexity. One of the results (Martin's theorem) was already mentioned. Even more general (in some sense) extension is discussed in~\cite{muchnik-game}.

\medskip
\textbullet\quad
It is easy to modify the proof to cover different versions of Kolmogorov complexity. For example, for prefix complexity we may  consider prefix-stable decompressors where $F(p)=x$ implies $F(p')=x$ for every $p'$ that has prefix $p$; similar modifications work for monotone and decision complexity. For \emph{a priori} complexity the players should specify lower approximations to a semimeasure.

\medskip
\textbullet\quad
One may change the rules of the game and let Alice and Bob directly provide upper bounds $KA$ and $KB$ instead of enumerating graphs for $A$ and $B$. Initially $KA(x)=KB(x)=+\infty$ for every $x$; at each step the player may decrease finitely many values of the corresponding function. The restriction (that goes back to Levin~\cite{levin-classification}) is that for every $n$ there is at most $2^n$ strings $x$ such that $KA(x)<n$ (the same restriction for $KB$). This approach works for prefix and decision complexities (but not for the monotone one).

\section{Examples}\label{sec:examples}

\subsection*{Conditional complexity and total programs}

Let $x$ and $y$ be two strings. The conditional complexity $\KS(x|y)$ of $x$ when $y$ is known can be defined as the length of the shortest program that transforms $y$ into $x$ (assuming the programming language is optimal). What if we require this program to be \textsl{total} (i.e., defined everywhere)?

It turns out that this requirement can change the situation drastically: there exist two strings $x$ and $y$ of length $n$ such that $\KS(x|y)=O(\log n)$ but any total program that transforms $y$ to $x$ has complexity (and length) $n-O(\log n)$. (Note that a total program that maps everything to $x$ has complexity at most $n+O(1)$, so the bound is quite tight.)

To prove this statement, we use the following game. Fix some $n$. We enumerate a graph of some function $f\colon\mathbb{B}^n\to\mathbb{B}^n$ (at each move we add some pairs to that graph). The opponent enumerates a list of at most $2^n-1$ total functions $g_1,g_2,\ldots$ (at each move the opponent may add some functions to this list). We win the game if there exist strings $x,y\in \mathbb{B}^n$ such that $f(y)=x$ but $g_i(y)\ne x$ for all $i$.

\textbf{Why we can win in this game}: First we choose some $x$ and $y$ and declare that $f(y)=x$. After every (non-trivial) move of the opponent we choose some $y$ where $f$ is still undefined and declare $f(y)=x$ where $x$ is different from currently known $g_1(y), g_2(y),\ldots$. The number of opponent's moves is less than $2^n$, therefore an unused $y$ still exists (we use only one point for every move of the opponent) and a value $x$ different from all $g_i(y)$ exists.

\textbf{Why the statement is true}: Let us use our strategy against the following opponent strategy: enumerate all total functions $\mathbb{B}^n\to\mathbb{B}^n$ that have complexity less than $n$. (Each function is considered here as a list of its values.) This strategy is computable (given $n$) and therefore the game is computable. Therefore, for the winning pair $(x,y)$ we have $\KS(x|y)=O(\log n)$ since $n$ is enough to describe the process and therefore to compute function $f$. On the other hand, any total function that maps $y$ to $x$ has complexity $n-O(1)$, otherwise the list of its values would appear in the enumeration.

So if we denote by $\overline{\KS}(x|y)$ the length of the shortest program for a total function that maps $y$ to $x$, we get a (non-computable) upper bound for $\KS(x|y)$ that sometimes differs significantly from $\KS$: it is possible that $\overline{\KS}(x|y)$ is about $n$  while $\KS$ is $O(\log n)$ (for strings $x$ and $y$ of length~$n$).

The conditional complexity defined is this way was considered also by Bruno~Bauwens~\cite{bruno-bauwens} (who used a different notation).

\subsection*{Extracting randomness requires $\Omega(\log n)$ additional bits}

Let us consider a question that can be considered as Kolmogorov-complexity version of
randomness extraction (though the similarity is superficial). Assume that a string $x$ is ``weakly random'' in the following sense: its complexity is high (at least $n$) but still can be much smaller than its length, which is polynomial in $n$. We want to ``extract'' randomness out of $x$, i.e., to get a string $y$ such that $y$ is random (=incompressible: its length is close to its complexity) using few additional bits, i.e., $\KS(y|x)$ should be small. When is it possible?

The natural approach: take the shortest program for $x$ as $y$. Then $y$ is indeed incompressible ($\KS(y)=l(y)+O(1)$; here $l(y)$ stands for the length of~$y$). And the complexity of $y$ when $x$ is known is $O(\log n)$: knowing $x$ and the length of a shortest program for $x$, we can find (at least some) shortest program for $x$. Taking the first $n$ bit of this shortest program, we get a string of length $n$, complexity $n+O(\log n)$ and  $O(\log n)$ conditional complexity relative to $x$.

What if we put a stronger requirement and requiere $\KS(y|x)$ to be $O(1)$ or $o(\log n)$? It turns that ``randomness extraction'' in this stronger sense is not always possible: \emph{there exists a string $x$ of length $n^2$ that has complexity at least $n$ such that every string $y$ of length $n$ that has conditional complexity $\KS(y|x)$ less than $0.5\log n$ has unconditional complexity $O(\log n)$} (i.e., is highly compressible). (The same result is true for all strings $y$ of length less than $n$, so we cannot extract even $n/2$ ``good random bits'' using $o(\log n)$ advice bits.)

To prove this statement, consider the following game. There are two sets $L=\mathbb{B}^n$ (``left part'') and $R=\mathbb{B}^{n^2}$ (``right part''). The opponent at each move may choose two elements $l\in L$ and $r\in R$ and add an edge between them (declaring $l$ to be a ``neighbor'' of $r$). The restriction is that every element in $R$ should have at most $d=\lceil \sqrt{n}\rceil$ neighbors. We may mark some elements of $L$ as ``simple''. We win if there is at least $2^n$ elements in $R$ that have the following property: \emph{all their neighbors are marked}.

\textbf{Why the statement is true} if we can win the game (using a computable strategy): Let the opponent declare $x\in L$ to be a neighbor of $y\in R$ if $\KS(x|y)<0.5\log n$. Then every $y$ has at most $d$ neighbors. The process is computable, so the game can be effectively simulated. Therefore, all $x$ declared as ``simple'' indeed have complexity $O(\log n)$ since each $x$ can be described by $n$ and its ordinal number in the enumeration of simple elements (the latter requires $0.5\log n$ bits). Among $2^n$ elements in $R$ that have the winning property there is one that has complexity at least $n$, and this is exactly what we claimed.

\textbf{How to win the game}: We do nothing while there are $2^n$ (or more) elements in $R$ that have no neighbors in $L$ (since this implies the required property). After $2^{n^2}-2^n$ elements get neighbors in $L$, we mark the neighbor that is used most often. It is a neighbor of at least $(2^{n^2}-2^n)/2^n=2^{n^2-n}-1 > 2^{n^2-2n}$ elements in $R$, and we restrict our attention to these ``selected'' elements ignoring all other elements of $R$. Then we do nothing while at least $2^n$ of  selected elements have no second neighbor. After that we mark the most used second neighbor and have at least $(2^{n^2-2n}-2^n)/2^n>2^{n^2-4n}$ elements that have two marked neighbors. In this way we either wait indefinitely at some step (and in this case we have at least $2^n$ elements that have only marked neighbors) or finally get $2^{n^2-2dn}>2^n$ elements who have $d$ marked neighbors and therefore cannot have non-marked ones, so we win.

Note that we could change the game allowing the opponent to declare $2^n$ elements in $R$ as simple and requiring in the winning condition that there is a non-simple element in $R$ that has no non-simple neighbors. This would make the game closer to original statement about Kolmogorov complexity but a bit more complicated.

This example is adapted from~\cite{ver-vyu}.

\subsection*{The compexity of a bijection}

For any two strings $x$ and $y$ one may look for a shortest program for a bijective function that maps $x$ to $y$. Evidently, it is not shorter than a shortest program for a total function that maps $x$ to $y$, therefore we get a lower bound $\overline{\KS}(y|x)-O(1)$ for a length (and complexity) of such a program. Since bijection can be effectively reversed, the bound can be made symmetric and we conclude that the length of a program for a bijection that maps $x$ to $y$ is at least
      $\max(\overline{\KS}(x|y), \overline{\KS}(y|x))-O(1).$
What about upper bounds?  Imagine there exists a simple total function that maps $x$ to $y$ and other simple total function that maps $y$ to $x$. Can we guarantee that there exists a simple \textsl{bijective} total function that maps $x$ to $y$?

To simplify the discussion, let us assume that $x$ and $y$ are of length $n$, the bijection should be length-preserving and $n$ is known (used as a condition in all the complexities).

This question corresponds to a game. Our opponent produces some total functions
     $$f_1,f_2,\ldots\colon\mathbb{B}^n\to\mathbb{B}^n \quad\hbox{and}\quad g_1,g_2,\ldots\colon\mathbb{B}^n\to\mathbb{B}^n$$
     claiming that one of $f_i$ maps $x$ to  $y$, and one of $g_j$ maps $y$ to $x$.  Knowing this functions (but not $x,y$), we have to produce bijections $$h_1,h_2,\ldots\colon\mathbb{B}^n\to\mathbb{B}^n$$ and guarantee that one of them maps $x$ to $y$.  (More precisely, the opponent wins if there exist $x$, $y$, $i$ and $j$ such that $f_i(x)=y$ and $g_j(y)=x$ but $h_k(x)\ne y$ for all $k$.) The question now is: how many bijections do we need to beat the opponent that can produce at most $m$ bijections of each type?

At first it seems that $m$ bijections are enough. Indeed, let us consider a bipartite graph where $x$ and $y$ are connected by an edge if  $f_i(x)=y$ and $g_j(y)=x$ for some $i$ and $j$. This graph has degree at most $m$ at both sides (e.g., $x$ can be connected only to $f_1(x),\ldots,f_m(x)$). Each bipartite graph where each vertex has degree at most $m$ and both parts are of the same size, can be covered by $m$ bijection graphs (we add edges to get degrees exactly $m$ and then use Hall's criterion for matchings).

This argument, if correct, would imply the upper bound $\max(\overline{\KS}(x|y),\overline{\KS}(y|x))+O(\log n)$ for the minimal complexity of the program that computes a bijection that maps $x$ to $y$. (Here $O(\log n)$ is added to take into account that we need to know $n$ for all our constructions.) Indeed, let the opponent to enumerate all the total functions $\mathbb{B}^n\to\mathbb{B}^n$ that have complexity at most
        $$u=\max(\overline{\KS}(x|y),\overline{\KS}(y|x)).$$
It is a computable process that involves at most $2^u$ functions. Beating this strategy of the opponent, we computably generate at most $2^u$ bijections (as we have assumed) and each bijection can be encoded by its ordinal number (at most $u$ bits) and $n$ (this requires $O(\log n)$ bits). Winning condition guarantees that one of these bijections maps $x$ to $y$.

However, this argument (and the result itself) is wrong. The problem is that the opponent does not tell us all its mappings at once but gives them one by one and we have to react immediately (otherwise we lose if the opponent does not make anything else). So we need to repeat this procedure after each move of the opponent, which gives $\Theta(m^2)$ bijection if opponent makes $m$ moves.

And this bound can be obtained by a much more simple strategy: for every $f_i$ and $g_j$ consider a bijection $h_{ij}$ that extents a partial matching
      $$x\leftrightarrow y \ \Leftrightarrow\ f_i(x)=y \ \text{and}\ g_j(y)=x.$$
This strategy gives upper bound $\overline{\KS}(x|y)+\overline{\KS}(y|x)+O(\log n)$.

\medskip

The main point of this example is that game arguments work in both directions: the absense of the winning strategy for us (and the existence of the winning strategy for the opponent) implies that the upper bound we wanted to prove is not true at all.

For example, the winning strategy in our game (for us) exists only if the number of our bijections is $\Omega(m^2)$ where $m$ is the maximal number of opponent's moves. It can be shown as follows. Let us assume that all the opponent's functions are constant functions (i.e., map all the elements of $\mathbb{B}^n$ into one element). In other terms, the opponent just selects vertices at both sides of the graph, and our goal is to provide bijections between each pair of selected vertices. It is easy to see that we would need $\Omega(m^2)$ bijections: indeed, if the opponent at each move selects a vertex that is not connected yet to  vertices selected earlier (which is always possible if the number of vertices is large compared to $m^2$) then we need $\Omega(m)$ new bijections to provide these new connections.

Translating this observation into Kolmogorov complexity language, we get the following statement: for every $k$ and $n$ such that $n>2k$ there exist two strings $x$ and $y$ of length $n$ such that $\KS(x),\KS(y)\le k+O(\log n)$ but any bijection that maps $x$ to $y$ has complexity $2k-O(1)$. To show this, use the trivial strategy at our side (we list all programs of length less than $2k$ that turn out to define a bijection $\mathbb{B}^n\to\mathbb{B}^n$; this property is enumerable) and let the opponent use the winning strategy described above (choosing elements not connected to already chosen elements by known bijections; the inequality $n>2k$ guarantees that $\Omega(2^k)$ steps are possible, since $(2^k)^2=2^{2k}< 2^n$). All chosen elements have complexity at most $k+O(\log n)$ and by the winning condition they are some of them not connected by a bijection of complexity less than $2k$.

\subsection*{Contrasting prefix and plain complexity}

Here we give a game-based proof of J.~Miller's result~\cite{miller-contrasting}. (The original proof in~\cite{miller-contrasting} uses a different scheme and involves the Kleene fixed-point theorem.)

Let $Q$ be a co-enumerable set of strings (i.e., its complement is enumerable) that for every $n$ contains at least one string of length $n$. Then \emph{for every $c$ there exists $n$ and $x$ of length $n$ such that $\KP(x)<n+\KP(n)-c$}. Here $\KP$ stands for prefix complexity; the contrast with the plain complexity arises because for plain complexity the set of incompressible strings (that have maximal possible complexity) is co-enumberable. (Note also that the maximal value of $\KP(x)$ for strings of length $n$ is $\KP(n)+O(1)$.)

To prove this statement, let us consider the following game specified by a natural number $C$ and a finite family of disjoint finite sets $S_1,\dots,S_N$. During the game each element  $s\in S=\cup_{j=1}^N S_j$ is labeled by two non-negative rational numbers $A(s)$ and $B(s)$ called ``Alice weight'' and ``Bob's weight''. Initially all weights are zeros. Alice and Bob make alternate moves. On each move each player may increase her/his weight of several elements~$s\in S$.

Both players must obey the following total weight restrictions:
       $$
\sum_{s\in S}A(s)\le1\quad\text{and}\quad \sum_{s\in S}B(s)\le1.
       $$
In addition, Bob must be ``fair'': for every $j$ Bob's weights of all $s\in S_j$ must be equal. That means that basically Bob assigns weights to $j\in\{1,\dots,N\}$ and Bob's weight $B(j)$ of $j$ is then evenly distributed among all $s\in S_j$ so that
      $$
B(s)=B(j)/\#S_j
      $$
for all $s\in S_j$.
Alice need not be fair.

This extra requirement is somehow compensated by allowing Bob to ``disable'' certain $s\in S$. Once an $s$ is disabled it cannot be ``enabled'' any more. Alice cannot disable or enable anything. For every $j$ Bob is not allowed to disable \textsl{all} $s\in S_j$:
every set $S_j$ should contain at least one element that is enabled (=not disabled).

The game is infinite. Alice wins if at the end of the game (or, better to say, in the limit) there exists an enabled $s\in S$ such that
      $$
\frac{A(s)}{B(s)}\ge C.
      $$

Now we have (as usual) to explain two things: why Alice has a (computable) winning strategy in the game (with some assumptions on the parameters of the game) and why this implies Miller's theorem.

\textbf{Lemma.} \emph{Alice has a computable winning strategy if $N\ge2^{8C}$ and $\#S_j\ge 8C$ for all $j\le N$}.

Let us show first why this statement implies the theorem. Let
     $$
C=2^{c} \quad\text{and}\quad N= 2^{8C}=2^{2^{c+3}}
    $$
Let us take the sets of all strings of length
           $$\log 8C+1,\dots, \log 8C+N$$
as $S_1,\ldots,S_N$. Then $S_j$ consists of $2^j\cdot 8C$ elements; the conditions of the lemma are satisfied and hence Alice has a computable winning strategy.

Consider the following Bob's strategy in this game: he enumerates the complement of $Q$ and disables all its elements; in parallel, he approximates the prefix complexity from above and once he finds out that $K(n)$ does not exceed some $l$, he increases the weights of all $2^n$ strings of length~$n$ up to $2^{-l-n}$. Thus at the end of the game $B(x)=2^{-K(n)-n}$ for all $s\in S$ that have length $n$ (i.e., for $s\in S_j$ where $j=n-\log 8C$).

Alice's limit weight function $x\mapsto A(x)$ is lower semi-computable given $c$, as both Alice's and Bob's strategies are computable given $c$. Therefore (since prefix complexity is equal to the logarithm of a priori probability)
   $$\KP(s|c)\le -\log A(s)+O(1)$$
for all $s\in S$. As Alice wins, there exists a string $s\in Q$ of some length $n\le N+\log 8C$ such that $A(s)/B(s)\ge C$, i.e.,
     $$
-\log A(s)\le -\log B(s)-c=\KP(n)+n-c.
    $$
This implies that
    $$
\KP(s|c)\le \KP(n)+n-c+O(1),
    $$
and
    $$\KP(s)\le \KP(n)+n-c+2\log c+O(1).$$
This is a bit weaker statement that we need: we wanted
    $$K(s)< K(n)+n-c.$$
To fix this, apply this argument to $c'=c+3\log c$ in place of $c$. For all large enough $c$ we then have $K(s)<K(n)+n-c$.

\smallskip

It remains to prove the Lemma by showing a winning strategy for Alice.

\smallskip

\textbf{Proof of the Lemma.}  The strategy is rather straighforward. The main idea is that playing with one $S_i$, Alice can force Bob to spend twice more weight than she does. Then she switches to next $S_i$ and so on until Bob's weight is exhausted while she has solid reserves. To achieve her goal on one set of $M$ elements, Alice assigns sequentially weights $1/2^M, 1/2^{M-1},\ldots,1/{2^1}$ and after each move waits until Bob increases his weight or disables the corresponding element. Since he cannot disable all elements and is forced to use the same weights for all elements while Alice puts more than half of the weight on the last element, Alice has factor $M/2$ as a handicap, and we may assume that $M$ beats $C$-factor that Bob has in his favor.

Now the formal details. Assume first that $\#S_j=M=4C$ for all $j$ and $N=2^{M}$. (We will show later how to adjust the proof to the case when $|S_j|\ge8C$ and $N\ge2^{8C}$.)

Alice picks an element $x_1\in S_1$ and assigns the weight $1/2^{M}$ to $x_1$.
Bob (to avoid losing the entire game) has either to assign a weight of more than $1/C2^{M}$ to all elements in $S_1$, or to disable $x_1$. In the second case Alice picks another element $x_2\in S_1$ and assigns a (twice bigger) weight of $2/2^{M}$ to it. Again Bob has a dilemma: either  to increase the weight for all elements of $S_1$  up to $2/C2^{M}$, or to disable $x_2$. In the second case Alice picks $x_3$, assigns a weight of $4/2^{M}$ to it, and so on. (If this process continues long enough, the last weight would be $2^{M-1}/2^M=1/2$.)

As Bob cannot disable all the elements of $S_1$, at some step $i$ the first case occurs, and Bob assigns a weight greater than $2^i/C2^M$ to all the elements of $S_1$.  Then Alice stops playing with $S_1$.  Note that the total Alice's weight of $S_1$ (let us call it $\beta$) is the sum of the geometric sequence:
        $$
\beta=1/2^{M}+2/2^M+\dots +2^{i-1}/2^M<2^i/2^M\le1.
       $$
Thus  Alice obeys the rules. Note that total Bob's weight of $S_1$ is more than $M2^{i-1}/C2^M=2^{i+1}/2^M$, which exceeds at least two times the total Alice's weight spent on $S_1$. This implies, in particular, that Bob cannot beat Alice's weight for the last element if the game comes to this stage (and Alice wins the game in this case.)

Then Alice proceeds to the second set $S_2$ and repeats the procedure. However this time she uses weights
     $
\alpha/2^{M},2\alpha/2^M,\dots,
    $
where $\alpha=1-\beta$ is the weight still available for Alice. Again she forces Bob to use twice more weight than she does. Then Alice repeats the procedure for the third set $S_3$ with the remaining weight etc.

Let  $\beta_j$  is the the total weight Alice spent on the sets $S_1,\dots,S_j$, and $\alpha_j=1-\beta_j$ the weight remaining after the first $j$ iterations. By construction, Bob's total weight spent on sets $S_1,\dots,S_j$ is greater than $2\beta_j$, so we have $2\beta_j<1$ and hence $\alpha_j> 1/2$. Consequently, Alice's total weight of each $S_j$ is more than $1/2^{M+1}$. Hence after at most $N=2^{M}$ iterations Alice wins.

If the size of $S_j$ are large but different, we need to make some modification. (We cannot use the same approach starting with $1/2^M$ where $M$ is the size of the set: if Bob beats the first element with factor $C$, he spends twice more weight than Alice but still a small amount, so we do not have enough sets for a contradiction.)

However, the modification is easy. If the number of elements in $S_j$ is a multiple of $4C$ (which is the case we use), we can split elements of $S_j$ into $4C$ groups of equal size, and treat all members of each group $G$ as one element. This means that if the above algorithm asks to assign to an ``element'' (group) $G$ a weight $w$, Alice distributes the weight $w$ uniformly among members of $G$ and waits until either Bob disables all elements of the group or assigns $4C$-bigger weight to all elements of $S_j$.

If $S_j$ is not a multiple of $4C$, the groups are not equal (the worst case is when some groups have one element while other have two elements), so to compensate for this we heed to use $8C$ instead of $4C$.

Note that excess in the number of groups (when $N$ is bigger than required $8C$) does not matter at all, we just ignore some groups.~\raisebox{-0.7ex}{\hbox{\small$\Box$}}

Note that this proof provides also some bound for $n$ (the length of the string); this bound is (almost) the same as given in Theorem~6.1 in~\cite{miller-contrasting}. Note also that instead of classifying strings according to their length, we could split them (effectively) into arbitrary finite sets $G_n$ whose cardinalities monotonically increase and are unbounded.  Then for every string $x\in G_n$ we have $\KP(x)\le \#G_n+KP(n)+O(1)$ and for every co-enumerable set $Q$ that intersects every $G_n$ there exists $n$ and $x\in G_N\cap Q$ such that $\KP(x)\le \#G_n+\KP(n)-c$ (for the same reasons).

\end{document}